\documentclass[12pt,twoside,reqno]{amsart}
\usepackage[colorlinks=true,citecolor=blue]{hyperref}
\usepackage{mathptmx, amsmath, amssymb, amsfonts, amsthm, mathptmx, enumerate, color,mathrsfs}
\usepackage{graphicx}

\usepackage{multirow}
\usepackage{epstopdf}
\usepackage{multicol}
\usepackage{algorithm}
\usepackage{algorithmic}
\usepackage{epstopdf}
\usepackage{subfig}

\newtheorem{theorem}{Theorem}[section]

\theoremstyle{definition}

\newtheorem{remark}[theorem]{Remark}
\numberwithin{equation}{section}

\input COTElvisMacros.mac
\begin{document}
\setcounter{page}{1}

%\vspace*{2.0cm}R
\title[Elvis Problem]
{A bisection method to solve the Elvis problem with convex bounded velocity sets}
\author[C. Graham, F. Marazzato \& P. Wolenski]{
Clinten A. Graham$^{1,2}$ \and Frederic Marazzato$^1$ \and Peter R. Wolenski$^1$ \\
\small{$^1$Department of Mathematics, Louisiana State University, Baton Rouge, LA 70803, USA}\\
   \small{email: \texttt{\{marazzato,pwolens\}@lsu.edu}\\} \\
  \small{$^2$Johns Hopkins University Applied Physics Laboratory, Laurel, MD, USA}\\
   \small{email: \texttt{clint.graham@jhuapl.edu}}}
\maketitle

\begin{abstract}
The Elvis problem has been studied in \cite{wolenski2021generalized}, which proves existence of solutions.
However, their computation in the non-smooth case remains unsolved.
A bisection method is proposed to solve the Elvis problem in two space dimensions for general convex bounded velocity sets.
The convergence rate is proved to be linear.
Finally, numerical tests are performed on smooth and non-smooth velocity sets demonstrating the robustness of the algorithm.
\end{abstract}

\section{Introduction}
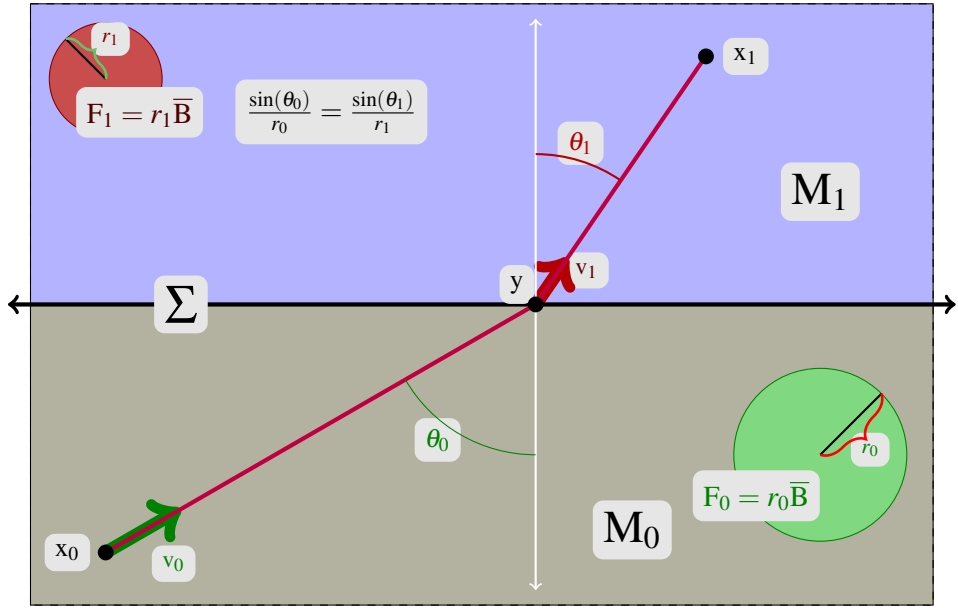
\begin{figure}[!htp]
\begin{center}
\begin{tikzpicture}[scale=1]

\coordinate (P1bdry) at (-6,4);
\coordinate (P2bdry) at (6,4);
\coordinate (P3bdry) at (-6,-4);
\coordinate (P4bdry) at (6,-4);

\coordinate (P1axis) at (-6,0);
\coordinate (P2axis) at (6,0);

\coordinate (P1) at (-5,-3.3);
\coordinate (P2) at (2.98,3.30);

\coordinate (Q) at (.716,0);
\coordinate (Q1) at (-1,0);

%Upper Zone
\draw[fill = blue!30]
(P1bdry) --  (P2bdry) -- (P2axis) -- (P1axis) -- cycle; 
%Lower Zone
\draw[fill = yellow!30!black!50]
(P1axis) --  (P2axis) -- (P4bdry) -- (P3bdry) -- cycle; 

%velocity balls
\fill[black!30!green!50] (4.5,-2) circle(1.15);
\draw[black!50!green] (4.5,-2) circle(1.15);
\draw[thick] (4.5,-2) -- +(45:1.15);
\node [black!50!green, fill=gray!20, rounded corners,below right] at (4.9,-1.7){\tiny{$r_0$}} ;
\draw[snake=brace,segment aspect=.5,segment amplitude=8,line width = 1pt,red] (4.5,-2)+(45:1.15) -- (4.5,-2);
\node [black!50!green, fill=gray!20, rounded corners,below left] at (-3.8,-3.2) {\scriptsize{$\v_0$}} ;
\node [black!50!green, fill=gray!20, rounded corners,below left] at (4.5,-2.2) {\small{$\F_0=r_0\oB$}} ;

\fill[black!30!red!70] (-5,3) circle(.75);
\draw[black!50!red] (-5,3) circle(.75);
\draw[thick] (-5,3) -- +(135:.75);
\node [black!50!red, fill=gray!20, rounded corners,above right] at (-5.2,3.3){\tiny{$r_1$}} ;
\draw[snake=brace,segment aspect=.5,segment amplitude=5,line width = 1pt,black!40!green!60!] (-5,3)+(135:.75) -- (-5,3);
\node [black!70!red, fill=gray!20, rounded corners,below right] at (1.1,.7){\scriptsize{$\v_1$}} ;
\node [black!70!red, fill=gray!20, rounded corners,below right] at (-5.4,2.9){\small{$\F_1=r_1\oB$}} ;

%light beam
\draw[->,black!50!green,line width=5pt] (P1) -- +(30:1.15);
\draw[->,black!30!red,line width=5pt] (Q) -- +(55.6:.75);
\draw[purple,ultra thick] (P1) -- (Q);
\draw[purple,ultra thick] (Q) -- ++(55.6:4);

%labels of manifolds
\node [fill=gray!20, rounded corners] at (2,-3){\Large{$\cM_0$}} ;
\node [fill=gray!20, rounded corners] at (4.5,1.5){\Large{$\cM_1$}} ;

%Boundaries
\draw[thin,dashed] (P1bdry)  -- (P2bdry) -- (P4bdry) -- (P3bdry) -- cycle;
\draw[<->,ultra thick] (-6.3,0) -- (6.3,0);

%angles
\draw[<->,white,thick] (.716,3.8) -- (.716,-3.8);
\node [black!30!red, fill=gray!20, rounded corners,above right] at (1,1.85)
{\footnotesize{$\theta_1$}};
\draw[black!30!red, thick] (.716,2) arc (90:55.6:2);
\node [black!50!green, fill=gray!20, rounded corners,above right] at (-.9,-2.1)
{\footnotesize{$\theta_0$}};
\draw[black!50!green, thin] (.716,-2) arc (270:210:2);

%points
\node [black, fill=gray!20, rounded corners,left] at (-5.2,-3.3) 
{\footnotesize{$\x_0$}};
\node [black, fill=gray!20, rounded corners,right] at (3.2,3.3) 
{\footnotesize{$\x_1$}};
\node [black, fill=gray!20, rounded corners,above left] at (.7,0){\footnotesize{$\y$}} ;
\fill[black] (P1) circle (3pt);
\fill[black] (Q) circle (3pt);
\fill[black] (P2) circle (3pt);

\node [fill=gray!20, rounded corners] at (-4,0){\LARGE{$\bSigma$}} ;
%Snell's
\node [black, fill=gray!20, rounded corners,below] at (-2,3) 
{\small{$\frac{\sin(\theta_0)}{r_0}=\frac{\sin(\theta_1)}{r_1}$}};
\end{tikzpicture}
\end{center}
\caption{The angles of incidence and Snell's Law}
\label{fig: Snell}
\end{figure}

Suppose $\cM_0,\,\cM_1\subseteq\R^2$ are, respectively, the lower and upper half-spaces in $\R^2$ and $\bSigma=\cl(\cM_0)\cap\cl(\cM_1)$ is the $x$--axis.
Fix $\mathbf{x}_0\in\cM_0$ and $\mathbf{x}_1\in\cM_1$. 
Suppose speed parameters $r_i>0$ are also fixed associated to the domain $\cM_i$ ($i=1,2$).
The setting is illustrated in Figure \ref{fig: Snell}.
The classical Snell's Law provides a necessary condition for a trajectory to traverse from $\x_0$ to $\x_1$ in least time while using maximal speed $r_i$ while in $\cM_0$.  The condition says 
\beq\label{eq: Snell}
\frac{\sin(\theta_0)}{r_0}=\frac{\sin(\theta_1)}{r_1},
\eeq 
where the $\theta_i$'s are the angles of incidence.
However, \eqref{eq: Snell} does not easily identify the point $\y\in\bSigma$ through which the optimal path passes. 
Rather, one usually sets this up as an elementary calculus problem as minimizing $x\mapsto \frac{|\mathbf{y}-\mathbf{x}_0|}{r_0}+\frac{|\mathbf{x}_1-\mathbf{y}|}{r_1}$ over $\mathbf{y}:=\bpm y\\ 0\epm\in\bSigma$.

This note provides an algorithm to find the optimal point $\mathbf{y}$ directly from \eqref{eq: Snell} but in a much more general situation that cannot in general be reduced to elementary calculus and instead relies on Convex Analysis.  

\section{Analysis of the problem}
We stay in $\R^2$, but similar results hold in any dimension $n$.  The centered balls as the (isotropic) velocity sets above are replaced by so-called {\em Elvis velocity} sets (which could be anisotropic) whose class is denoted by $\cC_0$.  A set $\F\in\cC_0$ is by definition nonempty, closed, convex, bounded, and contains $\0$ in its interior.  The least time a trajectory can go from $\0$ to some $\v\in\R^2$ using velocities from $\F$ is recorded by the gauge function
\[
\gamma_{\ \!\F}(\v):=\inf\left\{t> 0:\frac{1}{t}\v\in \F\right\}=\inf\left\{t> 0:\v\in t\F\right\},
\]
which is a finite-valued, positively homogeneous, convex function, see \cite{rockafellar2015convex}.  The generalized Elvis problem is
\begin{equation}\label{eq: gen CA}\tag*{$\bigl(P\bigr)$}
\inf\bigl[\gamma_{\ \!\F_0}(\mathbf{y}-\mathbf{x}_0)+\gamma_{\ \!\F_1}(\mathbf{x}_1-\mathbf{y})\bigr]\quad\text{over }\mathbf{y} \in\bSigma.
\end{equation}
This is a convex optimization problem.
Optimality conditions are contained in Theorem \ref{thm: Main} which is proved in \cite{wolenski2021generalized}.

\begin{thm}\label{thm: Main}
A necessary and sufficient condition for $\mathbf{y}\in\R^n$ to solve \ref{eq: gen CA} is the existence of $\mathbf{\zeta}_0$, $\bzeta_1\in\R^n$ satisfying
\begin{eqnarray}
\bzeta_0&\in &\partial\gamma_{\,\F_0}\bigl(\mathbf{y}-\mathbf{x}_0\bigr),\label{eq: polar0}  \\
-\bzeta_1&\in &\partial\gamma_{\,\F_1}\bigl(\mathbf{x}_1-\mathbf{y}\bigr),\quad\text{and} \label{eq: polar1}  \\
\bzeta_0+\bzeta_1 &\in & -\bN_{\bSigma}(\mathbf{y}).\label{eq: gen Snell}  
\end{eqnarray}
where $\bN_{\bSigma}(\mathbf{y})$ is the normal cone of $\bSigma$ at $\mathbf{y}$, which in this case is always the $y$-axis.
\end{thm}
\begin{remark}
Note that solutions $\mathbf{y} \in \mathbb{R}^n$ may not be unique if the velocity sets $F_0$ and $F_1$ are not strictly convex.
\end{remark}
We call \eqref{eq: gen Snell} the generalized Snell's Law from which the classical version
\eqref{eq: Snell} can be derived as follows.
Let us note the polar of $\F\in\cC_0$ as $\F^{\circ}$, which is defined by $\F^{\circ}:=\{\mathbf{\zeta} \in \mathbb{R}^n:\; \forall \mathbf{v} \in \F, \langle \bzeta,\mathbf{v}\rangle\leq 1\}$; it also belongs to $\cC_0$.
It is shown in \cite{wolenski2021generalized} that \eqref{eq: polar0} and \eqref{eq: polar1} imply that the optimal velocities are
\[ \mathbf{v}_0:=\frac{\mathbf{y}-\mathbf{x}_0}{\gamma_{\,\F_0}(\mathbf{y}-\mathbf{x}_0)}\in\F_0,
\quad \mathbf{v}_1:=\frac{\mathbf{x}_1-\mathbf{y}}{\gamma_{\,\F_1}(\mathbf{x}_1-\mathbf{y})}\in\F_1, \]
and that
\beq\label{eq: opt cond}
\bzeta_0\in\bN_{\,\F_0}\left(\mathbf{v}_0\right),\;
-\bzeta_1\in\bN_{\,\F_1}\left(\mathbf{v}_1\right),\;\text{and }
\gamma_{\,\F^{\circ}_0}(\bzeta_0)=1=\gamma_{\,\F^{\circ}_1}(-\bzeta_1).
\eeq
If $\F_0=r_0\oB$ and $\F_1=r_1\oB$, then $\F^{\circ}_0=\frac{1}{r_0}\oB$ and $\F^{\circ}_1=\frac{1}{r_1}\oB$.
The last statement in \eqref{eq: opt cond} then says 
\[
\bzeta_0=\frac{1}{r_0}\bpm \sin(\theta_0) \\ \cos(\theta_0)\epm \quad\text{and}\quad
-\bzeta_1=\frac{1}{r_1}\bpm \sin(\theta_1) \\ \cos(\theta_1)\epm 
\] 
for some angles $\theta_0$, $\theta_1$.
The condition \eqref{eq: gen Snell} says the $x$-component $\bzeta_0+\bzeta_1$ is $0$, which is \eqref{eq: Snell}.
The first condition in \eqref{eq: opt cond} implies $\u\mapsto\langle\bzeta_0,\u\rangle$ is maximized over $\u\in\F_0$ at $\mathbf{v}_0$ and the second condition in \eqref{eq: opt cond} implies $\u\mapsto\langle-\bzeta_1,\u\rangle$ is maximized over $\u\in\F_1$ at $\mathbf{v}_1$. 
The angles have therefore the same geometric meaning as in Figure~\ref{fig: Snell}.
However, the optimal $\mathbf{y}$ remains unknown.

\section{Bisection method}
We provide an algorithm to approximate solutions of \ref{eq: gen CA}.
For $\z=\bpm x \\ y\epm \in\R^2$, let $\Pi_\Sigma(\z)=x$ be the $x$-component of the orthogonal projection onto $\Sigma$.
It is assumed thereafter that $\Pi_\Sigma(\mathbf{x}_0) < \Pi_\Sigma (\mathbf{x}_1)$.
The goal of the algorithm is to compute the optimal $\mathbf{y} = \begin{pmatrix} y \\ 0 \end{pmatrix}$ that solves \ref{eq: gen CA}.
As the problem is not very regular, the existence of second order gradients is not assured.
Therefore, the gradient-free bisection method is chosen.

\subsection{The algorithm}
Let $\varepsilon>0$ be a tolerance parameter.
%Let $\phi = ...$ \TODO{Cite the paper of the ITP method. Also give the values of the parameters used.}
Algorithm \ref{al:bisection} is used to approximate a solution of \ref{eq: gen CA}.
\begin{algorithm}
\begin{enumerate}
\item Set $l^{0}:=\Pi_{\bSigma}(\mathbf{x}_0)$, $r^{0}:=\Pi_{\bSigma}(\mathbf{x}_1)$ and $y^0:=\frac12 \left(l^0+r^0 \right)$
\item  For $k \in \mathbb{N}^*$, let $\mathbf{y}^k = \begin{pmatrix} y^k \\ 0 \end{pmatrix}$ and
\begin{align*}
 \mathbf{v}_0^{k}:=\mathbf{y}^k -\mathbf{x}_0\quad&\text{and}\quad 
 \gamma_0^{k}:=\gamma_{\,\F_{0}}\left(\mathbf{v}_0^{k} \right)\\
 \mathbf{v}_1^{k}:=\mathbf{x}_1 - \mathbf{y}^k 
\quad &\text{and}\quad \gamma_1^{k}:=\gamma_{\,\F_{1}}\left(\mathbf{v}_1^{k} \right).
\end{align*}

\item Calculate (or choose any) $\bzeta_0^{k}$, $\bzeta_1^{k}$ so that
\begin{align*}
\bzeta_0^k \in \; \bN_{\F_0}\left({\mathbf{v}_0^{k}}\slash{\gamma_0^k}\right)\quad&\text{with}\quad\gamma_{\,\F_0^{\circ}}(\bzeta_0^k)=1 \\
-\bzeta_1^k \in \; \bN_{\F_1}\left({\mathbf{v}_1^{k}}\slash{\gamma_1^k}\right)\quad&\text{with}\quad\gamma_{\,\F_1^{\circ}}(-\bzeta_1^k)=1 
\end{align*}

\item Set $\delta=\Pi_{\bSigma}(\bzeta_0+\bzeta_1)$.  
\begin{enumerate}
\item If $|\delta| \le \varepsilon$, then the process terminates and $\mathbf{y}^k$ is the solution.
\item If $\delta < -\varepsilon$, then set $l^{k+1}:=l^{k}$, $y^{k+1}:=\frac12 \left(r^k + y^k \right)$ and $r^{k+1}:=y^{k}$, and start over at step 2. 
\item If $\delta > \varepsilon$, then set $r^{k+1}:=r^{k}$, $y^{k+1}:=\frac12 \left(l^k + y^k \right)$ and $l^{k+1}:=y^{k}$, and start over at step 2. 
\end{enumerate}
\end{enumerate}
\caption{Bisection}
\label{al:bisection}
\end{algorithm}

\begin{remark}
Note that if $\Pi_\Sigma (\mathbf{x}_0) > \Pi_\Sigma (\mathbf{x}_1)$, then the test $(b)$ and $(c)$ in step $(4)$ of the algorithm are reversed.
\end{remark}

\subsection{Convergence proof}
Let $d^k := r^k - l^k > 0$.
Thus, $d^0 = \Pi_{\bSigma}(\mathbf{x}_1)-\Pi_{\bSigma}(\mathbf{x}_0)$.
By construction, one has $d^{k+1} = \frac{d^k}2$.
Therefore,
\[ \bigl| \mathbf{y}^{k}-\mathbf{y}\bigr|
\leq d^{k} = \frac{d^0}{2^{k}} \]
Thus, $\mathbf{y}^{k} \to \mathbf{y}$ when $k\to\infty$ with a linear convergence rate.

\section{Numerical results}
We first verify that the convergence rate is linear and then provide an example with non-strictly convex and non-smooth velocity sets.
%The tolerance is set to $\varepsilon = 10^{-15}$, which is the machine error.

\subsection{Elliptic velocity sets}
%We take $\varepsilon$ to be the machine error.
%The test case consists in having $\mathbf{x}_0 := (-1,-1)$ and $\mathbf{x}_1 := (1,1)$.
%The velocity sets are the following ellipses $\F_0:=(\cos(t),\sin(t))$ and $\F_1:=(2\cos(t),\sin(t))$, for $t \in [0,2\pi]$.
%
%The green curve in Figure \ref{fig:ref ellipse}, which represents $\Pi_\Sigma (\bN_{\F_0} - \bN_{\F_1})$.
%The reference solution is computed as the abscissa of the unique point of the green curve in Figure \ref{fig:ref ellipse} which has a zero ordinate.
%\begin{figure}[!htp]
%\includegraphics[scale=0.5]{pics/ref_ellipse.pdf}
%\caption{Ellipcitc velocity sets: Reference solution.}
%\label{fig:ref ellipse}
%\end{figure}
%The reference solution is thus computed as $y^{\text{ref}} := -0.538264$.
%We write as $y^{\infty} := -0.538264$ the result of the algorithm and notice that $y^{\infty} = y^{\text{ref}}$.
%The convergence plot of the error is given in Figure \ref{fig:convergence ellipse}.
%\begin{figure}[!htp]
%\includegraphics[scale=0.5]{pics/convergence.pdf}
%\caption{Ellipcitc velocity sets: Convergence plot.}
%\label{fig:convergence ellipse}
%\end{figure}
%A linear regression gives a variance estimated at $5.4 \cdot 10^{-3}$ for a line of equation $\log(|y^k - y^\text{ref}|) = 0.0686556 - 0.00192272 k$ thus confirming the linear convergence of the algorithm.

% takin a more interesting set of ellipses
We take $\varepsilon$ to be the machine error.
The test case consists in having $\mathbf{x}_0 := (-1,-1)$ and $\mathbf{x}_1 := (1,1)$.
The velocity sets are the following ellipses $\F_0:=(\cos(t),\sin(t)/2)$ and $\F_1:=(2\cos(t),\sin(t))$, for $t \in [0,2\pi]$.

\begin{figure}[!htp]
	\includegraphics[scale=0.5]{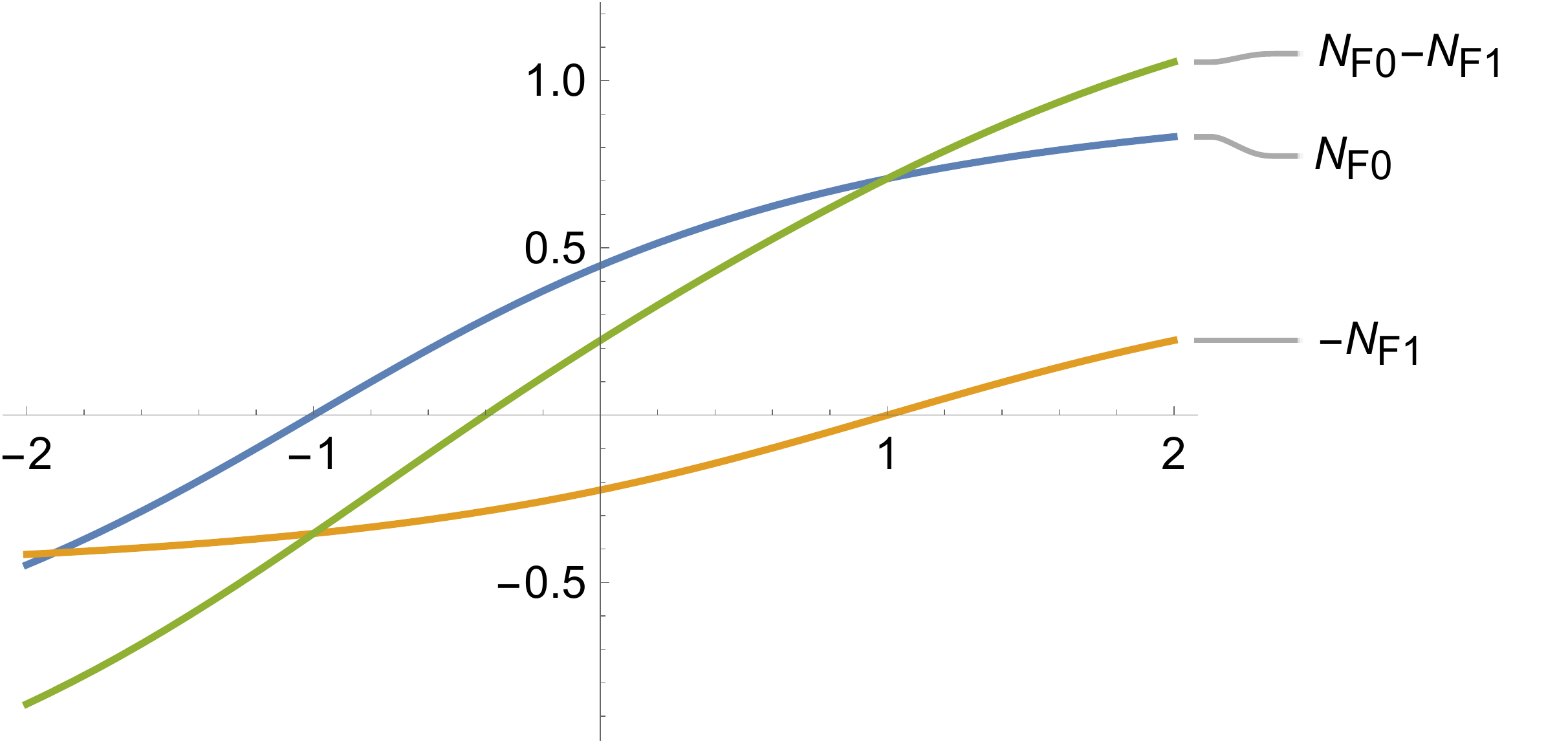}
	\caption{Ellipitc velocity sets: Reference solution.}
	\label{fig:ref ellipse}
\end{figure}

\begin{figure}[!htp]
	\includegraphics[scale=0.35]{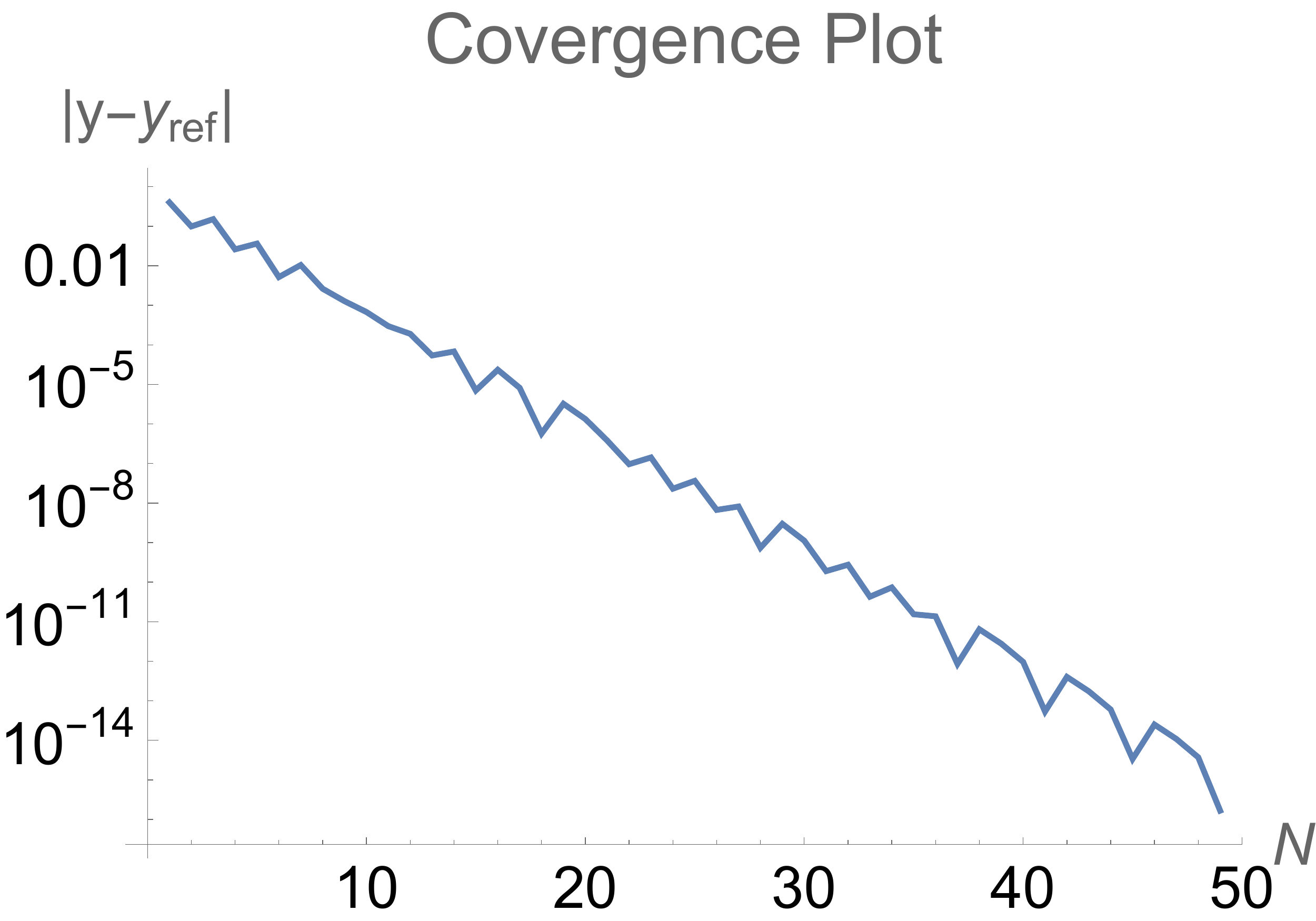}
	\caption{Ellipcitc velocity sets: Convergence plot.}
	\label{fig:convergence ellipse}
\end{figure}

The green curve in Figure \ref{fig:ref ellipse} represents $\Pi_\Sigma (\bN_{\F_0} - \bN_{\F_1})$.
The reference solution is computed as the abscissa of the unique point of the green curve in Figure \ref{fig:ref ellipse} which has a zero ordinate.
The reference solution is thus computed as %$y^{\text{ref}} := -0.401069$
 \[y^{\text{ref}} := -0.401. \]
The algorithm finds the result up to $\varepsilon$ as shown in the convergence plot of the error is given in Figure \ref{fig:convergence ellipse}.
A linear regression gives a variance estimated at $5.2 \cdot 10^{-3}$ for a line of equation $\log(|y^k - y^\text{ref}|) =-0.88 -0.0011 k$ thus confirming the linear convergence of the algorithm.

\subsection{Polyhedral velocity sets}
\begin{figure}[!htp]
	\includegraphics[scale=0.3]{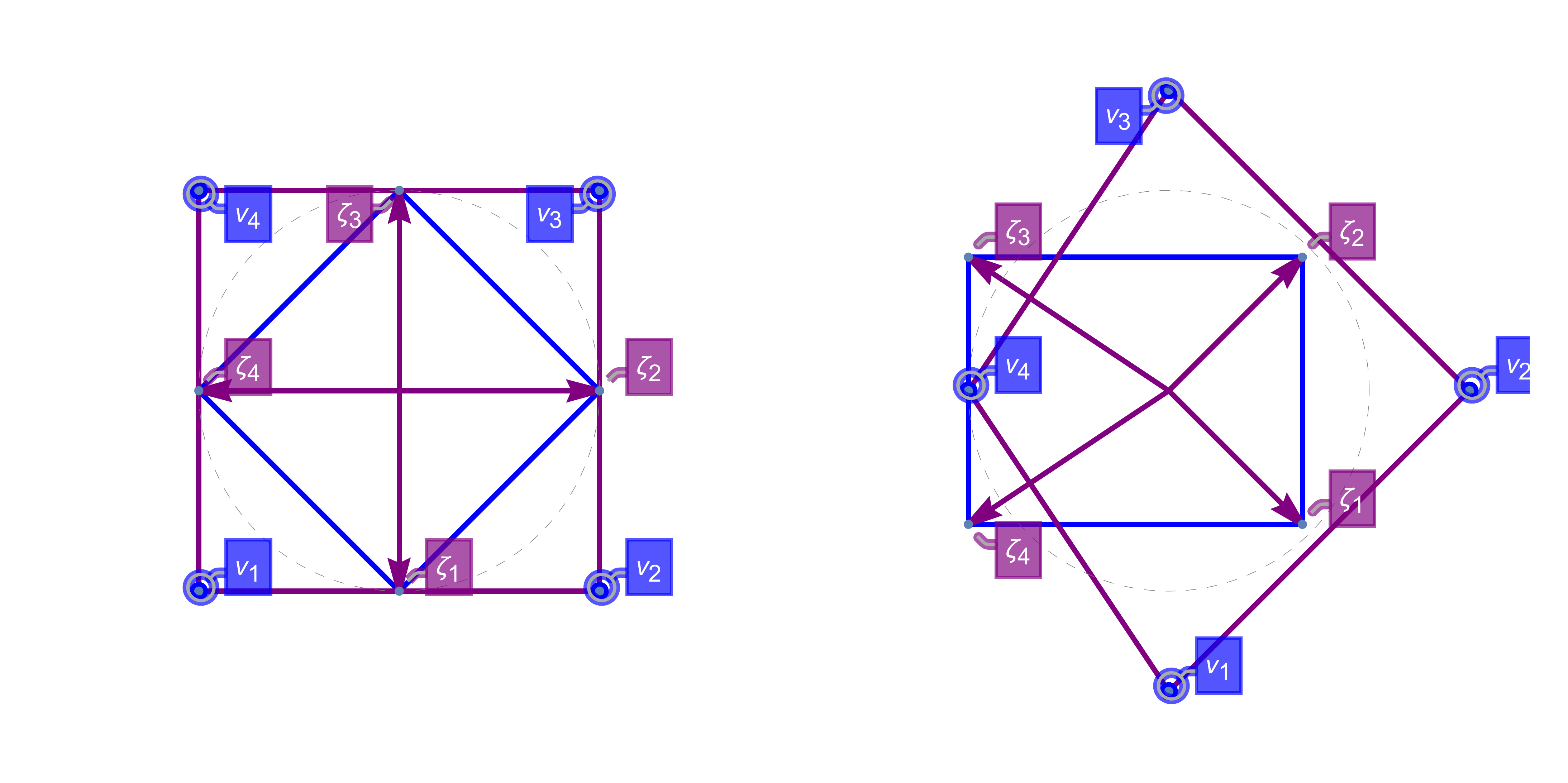}
	\caption{Polyhedral Velocity Sets: F0, F1 (left, right)}
	\label{fig:vel sets}
\end{figure}

\begin{figure}[!htp]
\includegraphics[scale=0.25]{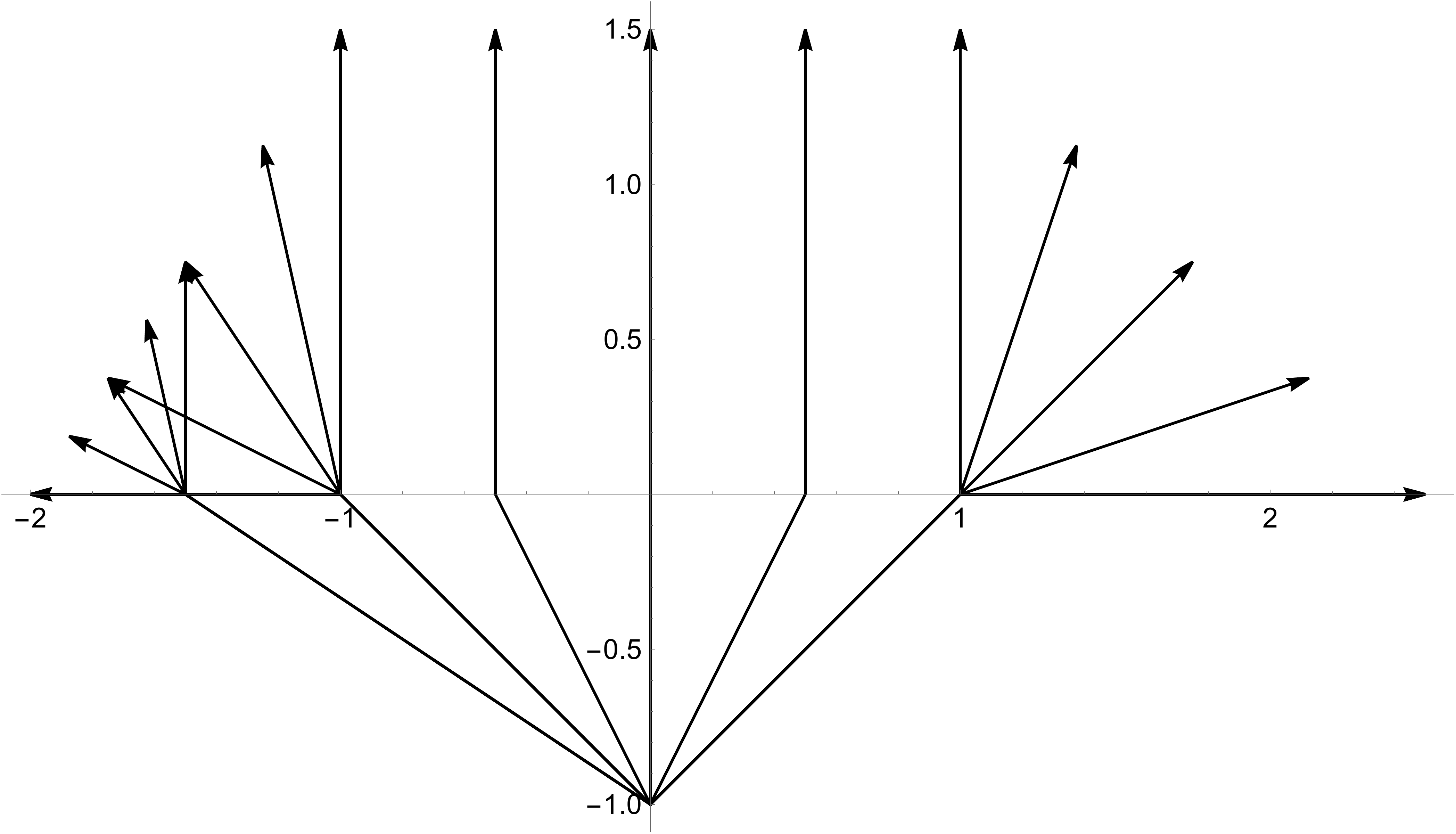}
\caption{Polyhedral Velocity Sets: Raycast Solutions}
\label{fig:raycast}
\end{figure}

%The test case consists in having $\mathbf{x}_0 := (0,-2)$ and $\mathbf{x}_1 := (6,10)$.
%The velocity sets are the squares plotted in Figure \ref{fig:vel sets}, both centred at the origin.
%\begin{figure}[!htp]
%\centering
%\subfloat{
%\begin{tikzpicture}[scale=1]
%\draw (1, 1) -- (-1,1) -- (-1,-1) -- (1, -1) -- cycle;
%\node at (0,0) {$\F_0$};
%\draw[<->] (-1.2, 1) -- (-1.2,-1);
%\node[left] at (-1.2, 0) {$1$};
%\end{tikzpicture}
%}
%\subfloat{
%\begin{tikzpicture}[scale=1]
%\draw (0, 1.5) -- (1.5,0) -- (0,-1.5) -- (-1.5, 0) -- cycle;
%\node at (0,0) {$\F_1$};
%\draw[<->] (0.2, 1.7) -- (1.7,0.2);
%\node[right] at (1, 1) {$1.5$};
%\end{tikzpicture}
%}
%\caption{Velocity sets.}
%\label{fig:vel sets}
%\end{figure}
%In that case, it can be shown that $\mathbf{y} = (0,2)$.
%\textcolor{red}{Give results.}
We take $\varepsilon$ to be the machine error. The test case consists in having $\mathbf{x}_0 := (0,-1)$ and $\mathbf{x}_1$ considered in several regions.
The velocity sets are the squares plotted in Figure \ref{fig:vel sets}, both containing the origin.These velocity sets were constructed to give rise to three regions of solutions, as can be seen in the raycast propagation \ref{fig:raycast}. %\newline

\begin{figure}[!htp]
	\includegraphics[scale=0.35]{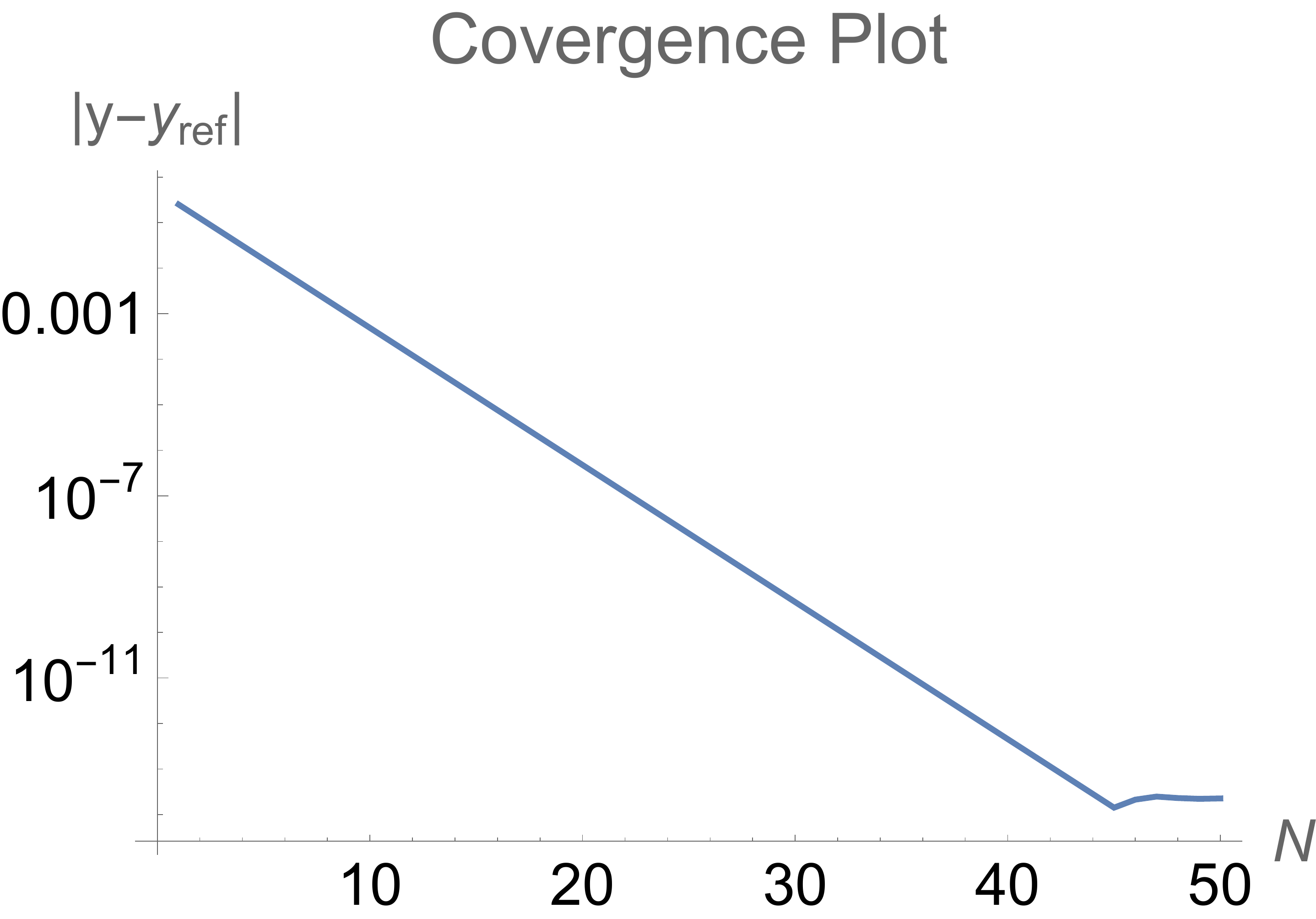}
	\caption{Polyhedral Velocity Sets: Convergence Plot}
	\label{fig:poly conv reg2}
\end{figure}

% region 1 (left)
When $\Pi_{\Sigma}(\mathbf{x}_1)<-1$ there is no unique solution. For $\mathbf{x}_1=(-2,1)$ the algorithm halts in two steps at $y^\text{ref}=-1.5$, the first solution it reaches in this region. %\newline
% region 2 (center)
When $-1\leq\Pi_{\Sigma}(\mathbf{x}_1)\leq1$ it can be shown that the optimal crossing point is $\Pi_{\Sigma}(\mathbf{x}_1)$. In this case the algorithm continuously bisects until it reaches $\varepsilon$ of the optimal solution. Figure \ref{fig:poly conv reg2} shows convergence when $\mathbf{x}_1=(0.5,1)$ and $y^\text{ref}=0.5$. %\linebreak

\section{Conclusion}
We have presented a bisection algorithm to solve the Elvis problem for general convex velocity sets.
Then we proved the convergence of the algorithm with a linear rate.
We then showed on the case of elliptical velocity sets that the proved convergence rate is optimal.
Finally, we applied the algorithm with non-smooth velocity sets.

Future work include generalizing this algorithm to $n$ space dimensions.
Also, a more efficient algorithm could be design on the model of the Newton--Raphson method to obtain quadratic convergence.
However, that will require the use of second order optimality conditions which are not common in convex analysis.

\section*{Code availability}
The algorithm has been implemented in Mathematica and is available at \url{https://github.com/clintg105/Elvis-Trajectory-Optimization.git}.

\bibliographystyle{plain}
\bibliography{bib}

\begin{thebibliography}{1}

\bibitem{rockafellar2015convex}
R.~T. Rockafellar.
\newblock {\em Convex analysis}.
\newblock Princeton university press, 2015.

\bibitem{wolenski2021generalized}
P.~R. Wolenski.
\newblock The generalized {E}lvis problem: Solving minimal time problems in
  anisotropic mediums.
\newblock In {\em 2021 60th IEEE Conference on Decision and Control (CDC)},
  pages 4552--4557. IEEE, 2021.

\end{thebibliography}

\end{document}